\documentclass{amsart}
\title{Pseudoholomorphic Curves in Nearly Kahler $\mathbf{CP}^3$}
\author{Feng Xu}
\newtheorem{theorem}{Theorem}[section]
\newtheorem{lemma}[theorem]{Lemma}

\theoremstyle{definition}

\newtheorem{corollary}[theorem]{Corollary}
\theoremstyle{remark}
\newtheorem{remark}[theorem]{Remark}

\numberwithin{equation}{section}

\begin{document}

\address{Department of Mathematics, Duke University, Durham, NC 27708}
\email{fxu@math.duke.edu}
\date{May 29, 2006}
\subjclass[2000]{Primary 53C38}
\keywords{Differential Geometry, Nearly Kahler, Pseudoholomorphic curves}
\begin{abstract}
We study pseudoholomorphic curves in the nearly Kalher $\mathbf{CP}^3$. It is shown that a class of curves called null-torsion are in one to one correspondence with the integrals of a holomorphic contact system on the usual Kahler $\mathbb{CP}^3$ studied by Bryant. Browing Bryant's result we get plenty of such curves. Rational curves are shown to be either vertical or horizontal or null-torsion. Both horizontal and null-torsion curves can be understood by mentioned Bryant's results while vertical curves are just fibers of the twistor fibration. 
\end{abstract}
\maketitle

\section{Introduction}
The most interesting almost Hermitian manifolds are perhaps the nearly Kahler manifolds. An almost Hermitian manifold is nealy Kahler provided its almost complex structure $J$ satisfies $\nabla_X(J)X=0$ for any vector field $X$ where $\nabla$ is the Levi-Civita connection. Although $J$ is nonintegrable, many aspects of Kahler geometry generalize to nearly Kahler manifold. For example, a generalized Hermitian-Yang-Mills theory was developed in \cite{Bryant3}. In this paper we are interested in another aspect, the pseudoholomorphic curves in nearly Kahler manifolds. Pseudoholomorphic curves in $S^6$ has been studied by Bryant in \cite{Bryant2}. He showed that every Riemann surface appears in $S^6$ as a null-torsion pseudoholomorphic curve with an arbitrarily large ramification degree. 

In this paper we study the next interesting case - curves in the  nearly Kahler $6$ manifold $\mathbf{CP}^3$. It is well known that there is a twistor fibration of the Kahler $3-$projective space, denoted by $\mathbb{CP}^3$, over $S^4$ with the fibers $\mathbb{CP}^1$. The nealy Kahler $\mathbf{CP}^3$ is defined by reversing the almost complex structure on the fibers. For each pseudoholomorphic curve $X: M^2\rightarrow \mathbf{CP}^3$ we construct holomorhpic sections $\mathbf{I}_1$, $\mathbf{I}_2$ and $\mathbf{II}$ of three holomorphic line bundles over $M$. We call those curves with $\mathbf{II}\equiv 0${\it null-torsion}. This is a well behaved condition. In fact, it turns out that null-torsion curves are in one to one correspondence with the holomorphic integrals of a holomorphic contact structure on the usual $\mathbb{CP}^3$. The latter has been thoroughly studied in \cite{Bryant1}. In particular, a Weierstrass formula was derived for such curves. Translating these results by this correspondence we get plenty of existence for null-torsion curves. If $M=S^2$ (these are called rational), it can be proved that $X$ necessarily falls into three categories: (1) vertical (2) horizontal (3) null-torsion. Both (2) and (3) can be reduced to the integrals of the holomorphic contact system mentioned before. Therefore we get a complete understanding of rational curves.

Fianlly, I would like to thank my advisor Robert Bryant for encouragements and stimulating discussions.       
\section{Structure Equations, Projective Spaces and the Flag Manifold} 
In this section we collect some facts needed in next section and formulate them in terms of the moving frame. Let $\mathbf{H}$ denote the real division algebra of quaternions. An element of $\mathbf{H}$ can be written uniquely as $q=z+jw$ where $z,w\in \mathbf{C}$ and $j\mathbf{H}$ satisfies 
  \[j^2=-1, zj=j\bar{z}\]
  for all $z\in\mathbf{C}.$ In this way we regard $\mathbf{C}$ as subalgebra of $\mathbf{H}$ and give $\mathbf{H}$ the structure of a complex vector space by letting $\mathbf{C}$ act on the right. 
  We let $\mathbf{H}^2$ denote the space of pairs $(q_1,q_2)$ where $q_i\in\mathbf{H}$. We will make $\mathbf{H}^2$ into a quarternion vector space by letting $\mathbf{H}$ act on the right 
  \[(q_1,q_2)q=(q_1q,q_2q).\]
This automatically makes $\mathbf{H}^2$ into a complex vector space of dimension $4$. In fact, regarding $\mathbf{C}^4$ as the space of $4-$tuples $(z_1,z_2,z_3,z_4)$ we make the explicit identification 
 \[(z_1,z_2,z_3,z_4)~(z_1+jz_2,z_3+jz_4).\] 
 This specific isomorphism is the one we will always mean when we write $\mathbf{C}^4=\mathbf{H}^2$. \\
 If $v\in\mathbf{H}^2\setminus(0,0)$ is given, let $v\mathbf{C}$ and $v\mathbf{H}$ denote respectively the complex line and the quarternion line spanned by $v$. The assigment $v\mathbf{C}\rightarrow v\mathbf{H}$ is a well defined mapping $T: \mathbf{CP}^3\rightarrow\mathbf{HP}^1$. The fibres of $T$ are $\mathbf{CP}^1$'s. So we have a fibration 
 \begin{equation}\label{fibration1}
 \begin{array}{lcc}
    \mathbf{CP}^1&\rightarrow &\mathbf{CP}^3\\
                 &            &\downarrow\\
                 &            &\mathbf{HP}^1\\ 
                 \end{array}\end{equation}
This is the famous twistor fibration. In order to study its geometry more thoroughly, we will now introduce the structure equations of $\mathbf{H}^2$.
First we endow $\mathbf{H}^2$ with a quarternion inner product $\left\langle, \right\rangle: \mathbf{H}^2\times\mathbf{H}^2\rightarrow \mathbf{H}$ defined by 
 \[\left\langle (q_1,q_2), (p_1,p_2)\right\rangle=\bar{q_1}p_1+\bar{q_2}p_2.\]
 We have identities
 \[\left\langle v,wq\right\rangle =\left\langle v,w\right\rangle q, \overline{\left\langle v,w\right\rangle}=\left\langle w,v\right\rangle,  \left\langle vq,w\right\rangle=\bar{q}\left\langle v,w\right\rangle.\]
 Moreover, $Re\left\langle,\right\rangle$ is a positive definite inner product which gives $\mathbf{H}^2$ the structure of a Euclidean space $\mathbf{E}^8$.
 Let $\mathfrak{F}$ denote the space of pairs $f=(e_1,e_2)$ with $e_i\in \mathbf{H}^2$ satisfying 
 \[\left\langle e_1,e_1 \right\rangle=\left\langle e_2, e_2\right\rangle=1, \left\langle e_1,e_2\right\rangle=0.\]
 We regard $e_i(f)$ as functions on $\mathfrak{F}$ with values in $\mathbf{H}^2$. Clearly $e_1(\mathfrak{F})=S^7\subset \mathbf{E}^8=\mathbf{H}^2$.
 It is well known that $\mathfrak{F}$ maybe canonically identified with $Sp(2)$ up to a left translation in $Sp(2)$. There are unique quaternion-valued 1-forms $\{\phi^{a}_{b}\}$ so that 
 \begin{equation}\label{strequation1}
 de_a=e_b\phi^b_{a},
 \end{equation}
 
  \begin{equation}\label{strequation2}
  d\phi^a_{b}+\phi^a_c\wedge\phi^c_b=0,
  \end{equation} 
    and 
    \begin{equation}
  \phi^a_{b}+\overline{\phi^b_a}=0.
  \end{equation}
  We have two canonical maps $C_1:\mathfrak{F}\rightarrow \mathbf{CP}^3$ and $C_2:\mathfrak{F}\rightarrow \mathbb{CP}^3$ by sending $f\in \mathfrak{F}$ to the complex lines spanned by $e_1(f)$ and $e_2(f)$ respectively. Recall that we have denoted the Kahler projective space by $\mathbb{CP}^3$ and the nearly Kahler one by $\mathbf{CP}^3$ whose structure will be explicitly described below. We are mainly interested in $\mathbf{CP}^3$. However, $\mathbb{CP}^3$ will play an important role. We now write structure equations for $C_1$ and $C_2$.
 First we immediately see that $C_1$ gives $\mathfrak{F}$ the structure of an $S^1\times S^3$ bundle over $\mathbf{CP}^3$ where  we have identified $S^1$ with the unit complex numbers and $S^3$ with the unit quaternions. The action is given by 
   \[f(z,q)=(e_1,e_2)(z,q)=(e_1z,e_2q),\]
 where $z\in S^1$ and $q\in S^3$. If we set
  $$\left[\begin{array}{ll} 
     \phi^1_{1}&\phi^1_2\\
     \phi^2_1 &\phi^2_2\\
     \end{array}\right]=
   \left[\begin{array}{cc}
     i\rho_1+j\overline{\omega_3}& -\frac{\overline{\omega_1}}{\sqrt{2}}+j\frac{\omega_2}{\sqrt{2}}\\
     \frac{\omega_1}{\sqrt{2}}+j\frac{\omega_2}{\sqrt{2}}& i\rho_2+j\tau\\
     \end{array}\right]   
   $$
  where $\rho_1$ and $\rho_2$ are real 1-forms while $\omega_1$, $\omega_2$, $\omega_3$ and $\tau$ are complex valued,  we may rewrite one part of the structure equation (\ref{strequation2}) relative to the $S^1\times S^3$ structure on $\mathbf{CP}^3$ as
   \begin{equation}\label{strequation4}d\left(
   \begin{array}{l}
   \omega_1\\
    \omega_2\\
    \omega_3\\
    \end{array}\right)=-\left(
    \begin{array}{lcl}
    i(\rho_2-\rho_1)&-\bar{\tau}&0\\
    \tau& -i(\rho_1+\rho_2)&0\\
    0&0&2i\rho_1\\
    \end{array} \right)\wedge\left(
    \begin{array}{l}
    \omega_1\\\omega_2\\\omega_3\end{array}\right)+\left(
    \begin{array}{l}    \overline{\omega_2\wedge\omega_3}\\\overline{\omega_3\wedge\omega_1}\\\overline{\omega_1\wedge\omega_2}\end{array}\right).
 \end{equation}  
 This in particular defines a nearly Kahler structure on $\mathbf{CP}^3$ by requiring $\omega_1$, $\omega_2$ and $\omega_3$ to be of type $(1,0)$ (note that this almost complex structure is nonintegrable, thus different from the usual integrable one). We denote
 $$\left(\begin{array}{ll}
    \kappa_{1\overline{1}}&\kappa_{1\bar{2}}\\
    \kappa_{2\bar{1}}&\kappa_{2\bar{2}}\\
    \end{array}\right)=\left(
    \begin{array}{cl}
    i(\rho_2-\rho_1)&-\bar{\tau}\\
    \tau&-i(\rho_1+\rho_2)\\
    \end{array}\right)
    $$
    and $\kappa_{3\bar{3}}=2i\rho_1$ in the usual notation of a connection. Then the other part of the structure equation (\ref{strequation2}) may be written as the curvature of this nearly Kahler structure
    \begin{equation}\label{strequation5}d\left(
    \begin{array}{lc}
    \kappa_{1\bar{1}}&\kappa_{1\bar{2}}\\
    \kappa_{2\bar{1}}&\kappa_{2\bar{2}}\\
    \end{array}\right)+\left(\begin{array}{lc}
    \kappa_{1\bar{1}}&\kappa_{1\bar{2}}\\
    \kappa_{2\bar{1}}&\kappa_{2\bar{2}}\\
    \end{array}\right)\wedge\left(\begin{array}{lc}
    \kappa_{1\bar{1}}&\kappa_{1\bar{2}}\\
    \kappa_{2\bar{1}}&\kappa_{2\bar{2}}\\
    \end{array}\right)
    =
    \left(\begin{array}{lc}
    \omega_1\wedge\overline{\omega_1}-\omega_3\wedge\overline{\omega_3}& \omega_1\wedge\overline{\omega_2}\\
    \omega_2\wedge\overline{\omega_1}&\omega_2\wedge\overline{\omega_2}-\omega_3\wedge\overline{\omega_3}\\
    \end{array}\right),
    \end{equation}
 as well as
 \begin{equation}\label{strequation6}
 d\kappa_{3\bar{3}}=-(\omega_1\wedge\overline{\omega_1}+\omega_2\wedge\overline{\omega_2}-2\omega_3\wedge\overline{\omega_3}).
 \end{equation}
 
In an exactly analogous fashion, $C_2$ gives $\mathfrak{F}$ a structure of an $S^1\times S^3$ bundle over $\mathbb{CP}^3$ with the action now given by
   \[(e_1,e_2)(q,z)=(e_1q,e_2z),\]
where $z\in S^1$ and $q\in S^3$. However, $\omega_1$, $\omega_2$, $\kappa_{1\bar{2}}$ and their complex conjugates become semibasic and $\omega_3$ is not. The usual Kahler structure on $\mathbb{CP}^3$ is defined by requiring $\frac{\overline{\omega_1}}{\sqrt{2}}$, $\frac{\omega_2}{\sqrt{2}}$ and $\kappa_{2\bar{1}}$ to be of type $(1,0)$ and unitary. Relative to this Kahler structure, we may rewrite part of the structure equations as
\begin{equation}
d\left(
\begin{array}{l}
\frac{\overline{\omega_1}}{\sqrt{2}}\\
\frac{\omega_2}{\sqrt{2}}\\
\kappa_{2\bar{1}}
\end{array}\right)=-\left(
\begin{array}{lcc}
-\kappa_{1\bar{1}}&\omega_3&-\frac{\overline{\omega_2}}{\sqrt{2}}\\
-\overline{\omega_3}&\kappa_{2\bar{2}}&\frac{\omega_1}{\sqrt{2}}\\
\frac{\omega_2}{\sqrt{2}}&-\frac{\overline{\omega_1}}{\sqrt{2}}&\kappa_{2\bar{2}}-\kappa_{1\bar{1}}\\

\end{array}\right)\left(\begin{array}{l}
\frac{\overline{\omega_1}}{\sqrt{2}}\\
\frac{\omega_2}{\sqrt{2}}\\
\kappa_{2\bar{1}}
\end{array}\right).
\end{equation}
 We will also need some properties of the flag manifold $\mathbf{Fl}=\mathfrak{F}/(U(1)\times U(1))$. Equivalently $\mathbf{Fl}$ consists of pairs of complex lines $([e_1], [e_2])$ with $\left\langle e_1, e_2\right\rangle=0$. Of course $\mathfrak{F}$ defines a natural $S^1\times S^1$ structure on $\mathbf{Fl}$ for which the forms $\omega_1,\omega_2,\omega_3,\kappa_{2\bar{1}}$ and their complex conjugates are semibasic. Moreover, we have a double fibration of $\mathbf{Fl}$ over the two projective spaces:
\begin{equation}\label{doublefibration}
\begin{array}{lllll}
 & &\mathbf{\mathfrak{F}}& & \\
 & &\downarrow& &\\
 & &\mathbf{Fl}& &\\
 &\swarrow&&\searrow&\\
\mathbf{CP}^3&  & &&\mathbb{CP}^3\\
\end{array}\end{equation}
We denote the first fibration by $\Pi_1$ and the second fibration by $\Pi_2$. Explicitly $\Pi_a$ ($a=1, 2$) sends $([e_1], [e_2])\in\mathbf{Fl}$ to the complex line $[e_a]$. By requiring $\overline{\omega_1},\omega_2,\omega_3,\kappa_{2\bar{1}}$ to be complex linear we define an almost complex structure on $\mathbf{Fl}$. It is easy to check from the structure equations that this almost complex structure is integrable and $\Pi_2$ is thus a holomorphic projection.

Finally there are various complex vector bundles associated with $\mathfrak{F}$ which will be important. First on $\mathbf{CP}^3$, there are two obvious complex bundles, the tautological bundle $\epsilon$ and the trivial rank $4$ bundle $\mathbf{C}^4$. We view $\epsilon$ as the subbundle of $\mathbb{C}^4$ spanned by $e_1$ and denote the quotient bundle by $Q$. Using the obvious Hermitian product, we identify $Q$ as a subbundle of $\mathbf{C}^4$ locally spanned by $e_1j$, $e_2$ and $e_2j$. Note that $e_1j$ itself spans a well-defined line bundle, which is isomorphic to $\epsilon^*$. Denote the quotient $Q/\epsilon^*$ by $\tilde{Q}$ which, again, may be regarded as a subbunle of $\mathbf{C}^4$ locally spanned by $e_2, e_2j$.  We write $T\mathbf{CP}^3$ the complex tangent bundle of $\mathbf{CP}^3$. The $S^1\times S^3$ determines a splitting 
   \[T\mathbf{CP}^3=\mathcal{H}\oplus \mathcal{V},\]
   where $\mathcal{H}$ has rank $2$ and $\mathcal{V}$ has rank $1$. We call $\mathcal{H}$ the horizontal part and $\mathcal{V}$ the vertical part relative to the fibration \ref{fibration1}. One can show $\mathcal{V}$ is isomorphic to $\epsilon^2$ as a Hermitian line bundle by locally identifying $\frac{1}{\sqrt{2}}e_1\otimes e_1$ with the complex tangent vector dual to the $(1,0)$ form $\omega_1$, denoted by $f_{\bar{1}}$. Similarly $\mathcal{H}$ is isomorphic to $\epsilon^*\otimes \tilde{Q}$ with $\frac{1}{\sqrt{2}}e_1^*\otimes e_2$ identified with the tangent vector $f_{\bar{2}}$ dual to $\omega_2$ and $\frac{1}{\sqrt{2}}e_1^*\otimes e_2j$ identified with $f_{\bar{3}}$ dual to $\omega_3$. Pulled back to $\mathbf{Fl}$ by $\Pi_1$, the bundles $\tilde{Q}$ and $\mathcal{H}$ splits as 
    \[\Pi_1^*\tilde{Q}=\tilde{\epsilon}\oplus\tilde{\epsilon}^*,\]
    where $\tilde{\epsilon}$ is locally spanned by $e_2$ and $\tilde{\epsilon}^*$ denotes its dual, locally spanned by $e_2j$. Of course $\Pi_1^*\mathcal{H}$ splits correspondingly as 
    \[\Pi_1^*\mathcal{H}=\epsilon^*\otimes\tilde{\epsilon}\oplus\epsilon^*\otimes\tilde{\epsilon}^*.\]
Similar constructions apply to $\mathbb{CP}^3$. We will only point out differences and some relations. The tautological bundle on $\mathbb{CP}^3$ becomes $\tilde{\epsilon}$ when pulled back to $\mathbf{Fl}$. The complex tangent bundle also splits as a sum of a vertical part $\mathbb{V}$ and 
a horizontal part $\mathbb{H}$. The vertical part is isomorphic to $(\tilde{\epsilon}^*)^2$ compared with $\mathbf{CP}^3$ case because of the reversed almost complex structure. The horizontal part, when pulled back by $\Pi_2$ to $\mathbf{Fl}$ is isomorphic to $\tilde{\epsilon}^*\otimes\epsilon^*\oplus\tilde{\epsilon}^*\otimes\epsilon$. Note this splitting shares a common factor with $\Pi_1^*\mathcal{H}$ which will become important later. 
In the various isomorphisms, we no longer need the $\frac{1}{\sqrt{2}}$ to make them Hermitian.
Moreover, since the comples structure on $\mathbb{CP}^3$ is integrable, many of these bundles have holomorphic structures. Among them the dual of the vertical tangent bundle of $\mathbb{CP}^3$, which we denote by $\mathbb{V}^*$ is particularly important. Locally $\mathbb{V}^*$ is spanned by $\kappa_{2\bar{1}}$ as a subbundle of the complex cotangent bundle of $\mathbb{CP}^3$. We have the following result due to R. Bryant \cite{Bryant1}
\begin{lemma}
The bundle $\mathbb{V}^*$ is isomorphic to $\tilde{\epsilon}^2$ as a Hermitian holomorphic line bundle. Moreover, it induces a holomorphic contact structure on $\mathbb{CP}^3$.
\end{lemma} 
 The integrals of this holomorphic contact system was thoroughly investigated in \cite{Bryant1} (see Section 3).

 \section{Pseudoholomorphic Curves in $\mathbf{CP}^3$}
 Let $M^2$ be a connected Riemann surface. A map $X:M^2\rightarrow \mathbf{CP}^3$ is called a pseudoholomorphic curve if $X$ is nonconstant and the differential of $X$ commutes with the almost compex structures.
 We let $x:\mathfrak{F}_X\rightarrow M^2$, $\mathcal{V}_X\rightarrow M^2$ and $\mathcal{H}_X\rightarrow M^2$ be the pull back bundles of $\mathfrak{F}\rightarrow \mathbf{CP}^3$, $\mathcal{V}\rightarrow \mathbf{CP}^3$ and $\mathcal{H}\rightarrow \mathbf{CP}^3$ respectively. Thus for intance, we have
   \[\mathfrak{F}_X=\{(x,f)\in M^2\times \mathfrak{F}|X(x)=C_1(f)\}.\]
   Of course, $\mathfrak{F}_X$ is an $S^1\times S^3$ bundle over $M^2$ and $\mathcal{V}_X$ and $\mathcal{H}_X$ are Hermitian complex bundles of rank $1$ and $2$ respectively. Moreover, the natural map $\mathfrak{F}_X\rightarrow \mathfrak{F}$ pulls back various quatities on $\mathfrak{F}$, which we still denote by the same letters. For example, $f_{\bar{1}}, f_{\bar{2}}$ now denote functions on $\mathfrak{F}_X$ valued in $\mathcal{H}_X$. The structure equations (\ref{strequation4}), (\ref{strequation5}) and (\ref{strequation6}) still hold, on $\mathfrak{F}_X$ now. Also for functions and sections with domains in $M^2$, we will pull these back up via $x^*$ to $\mathfrak{F}_X$. For example, any section $s:M^2\rightarrow \mathcal{H}_X$ can be written in the form $s=f_{\bar{1}}s_1+f_{\bar{2}}s_2$ where $s_i$ are complex functions on $\mathfrak{F}_X$. Using this convention, the pullback of $\kappa$ induces  connections on $\mathcal{H}_X$ and $\mathcal{V}_X$ compatible with the Hermitian structures. Namely $\nabla: \Gamma(\mathcal{H}_X)\rightarrow \Gamma(\mathcal{H}_X\otimes T^*M^2)$ is given by 
    \[\nabla(f_{\bar{i}}s_i)=f_{\bar{i}}\otimes(ds_i+\kappa_{i\bar{j}}s_j).\]
Since we are working over a Riemann surface, it is well-known that there are unique holomorphic structures on $\mathcal{H}_X$ and $\mathcal{V}_X$ compatible with these connections. From now on we will regard these two bundles as holomorphic Hermitian vector bundles over $M^2$.\\
Another thing to notice is that $\{\omega_i\}$ are semi-basic with respect to $x:\mathfrak{F}_X\rightarrow M^2$. Moreover, they are of type $(1,0)$ since $dX$ is complex linear. Set \[\mathbf{I}_1=f_{\bar{1}}\otimes\omega_1+f_{\bar{2}}\otimes\omega_2, \mathbf{I}_2=f_{\bar{3}}\otimes\omega_3.\] It is clear that $\mathbf{I}_1$ and $\mathbf{I}_2$ are well defined sections of $\mathcal{H}_X\otimes T^*M^2$ and $\mathcal{V}\otimes T^*M^2$ respectively where $T^*M^2$ is the holomorpic line bundle of $(1,0)$ forms on $M^2$.
 \begin{lemma}\label{LemmaI}
 The sections $\mathbf{I}_1$ and $\mathbf{I}_2$ are holomorphic. Moreover, $\mathbf{I}_1$ and $\mathbf{I}_2$ only vanish at isolated points unless $X(M^2)$ is horizontal ( when $\mathbf{I}_2$ vanishes identically) or vertical (when $\mathbf{I}_1$ vanishes identically and thus $X(M^2)$ is an open set of a fiber $\mathbf{CP}^1$ in (\ref{fibration1})).
 \end{lemma}
 \begin{proof}
 We only show $\mathbf{I}_1$ is holomorphic and leave $\mathbf{I}_2$ for the reader. Choose a uniformizing parameter $z$ on a neighborhood of $x_0\in M$. In a neighborhood of $x^{-1}(x_0)$, there exist functions $a_i$ so that $\omega_i=a_idz$. It follows that $\omega_i\wedge\omega_j=0$, so we have $d\omega_i=-\kappa_{i\bar{j}}\wedge\omega_j$. This translates to $(da_i+\kappa_{i\bar{j}}a_j)\wedge dz=0$ so there exists $b_i$ so that 
  \[da_i+\kappa_{i\bar{j}}a_j=b_idz.\]
Thus, when we compute $\overline{\partial}\mathbf{I}_1$ we have 
$$\begin{array}{lcl}
\bar{\partial}\mathbf{I}_1&=&(\nabla(f_{\overline{i}}a_i)\otimes dz)^{0,1}\\
                 &=&f_{\bar{i}}\otimes dz\otimes (da_i+\kappa_{i\bar{j}}a_j)^{0,1}\\
                 &=&f_{\bar{i}}\otimes dz\otimes (b_idz)^{0,1}\\
                 &=&0,
\end{array}$$
 so $\mathbf{I}_1$ is holomorphic.
 Moreover, by complex analysis, if $\mathbf{I}_1$ or $\mathbf{I}_2$ vanishes at a sequence of points with an accumulation, the section has to be identically $0$ since $M^2$ is connected. 
 
\end{proof} 
\begin{remark}
It is clear that $\mathbf{I}_1$ and $\mathbf{I}_2$ are just horizontal and vertical parts of the evaluation map 
$X_*(TM)\rightarrow T_X$.
\end{remark}
We will call a curve with $\mathbf{I}_1=0$ ($\mathbf{I}_2=0$) vertical (horizontal). Of course vertical curves are just the fibers $\mathbf{CP}^1$ of $T$. To study horizontal curves it does no harm to reverse the almost complex structure on the fiber of $T$. This new complex structure is integrable and actually equivalent to the usual complex structure on the $3$ projective space. 
The horizontal bundle $\mathcal{H}$ turns out to be a holomorphic contact structure under the usual complex structure. The integral curves of this contact system are thoroughly described in \cite{Bryant1}. We therefore have a good understanding of horizontal pseudoholomorphic curves in $\mathbf{CP}^3$.

We now assume both $\mathbf{I}_1$ is not identically $0$. There exists a holomorphic line bundle $L\subset \mathcal{H}_X$ so that $\mathbf{I}_1$ is a nonzero section of $L\otimes T^*M$. We let $R_1$ be the ramification divisor of $\mathbf{I}_1$. That is, 
\[ R_1=\sum_{p: \mathbf{I}_1(p)=0}ord_p(\mathbf{I}_1)p.\]
 $R_1$ is obviously effective, and we have 
 \[L=TM\otimes [R_1].\] 
 Similarly if $\mathbf{I}_2$ does not vanish identically let $R_2$ be the ramification divisor of $\mathbf{I}_2$.  Then $R_2$ is effective and 
 \[\mathcal{V}_X=TM\otimes [R_2].\]
 Now we adapt frames in accordance with the general theory. We let $\mathfrak{F}^{(1)}_X$ be the subbundle of pairs $(x,f)$ with $f_{\bar{2}}\in L_x$. Then $\mathfrak{F}^{(1)}_X$ is a $U(1)\times U(1)$ bundle over $M$. The cononical connection on $L$ is described as follows: If $s:M\rightarrow L$ is a section, then $s=f_{\bar{2}}s_2$ for some function $s_2$ on $\mathfrak{F}^{(1)}_X$. Then
    \[\nabla s=f_{\bar{2}}\otimes (ds_2+\kappa_{2\bar{2}}s_2).\]
 Similarly the quotient bundle $N_X=\mathcal{H}_X/L$ has a natural holomorphic Hermitian structure. Let $(f_{\bar{1}}):\mathfrak{F}^{(1)}_X\rightarrow N_X$ be the function $f_{\bar{1}}$ followed by the projection $\mathcal{H}_X\rightarrow N_X$. If $s:M\rightarrow N_X$ is any section, then $s=(f_{\bar{1}})s_1$ for $s_1$ on $\mathfrak{F}^{(1)}_X$ and we have 
   \[\nabla s=(f_{\bar{1}})\otimes (ds_1+\kappa_{1\bar{1}}s_1).\]
   Note since $\mathbf{I}_1$ has values in $L\otimes T^*M$, we must have $\omega_1=0$ on $\mathfrak{F}^{(1)}_X$. If we differential this using structures (\ref{strequation4}) we have 
    \[d\omega_1=-\kappa_{1\overline{2}}\wedge\omega_2=0.\]
 It follows that $\kappa_{1\bar{2}}$ is of type $(1,0)$.
 \begin{lemma}
 Let $\mathbf{II}=(f_{\bar{1}})\otimes f_{2}\otimes\kappa_{1\bar{2}}$ where $f_2$ is the dual of $f_{\bar{2}}$. Then $\mathbf{II}$ is a holomorphic section of $N_X\otimes L^*\otimes T^*M$.
  \end{lemma}
  \begin{proof}
  Since $\kappa_{1\bar{2}}$ is of type $(1,0)$, there exists $b$ locally such that $\kappa_{1\overline{2}}=bdz$. The structure equations (\ref{strequation5}) pulled back to $\mathfrak{F}^{(1)}$ gives
   $d\kappa_{1\bar{2}}=-\kappa_{1\bar{1}}\wedge\kappa_{1\bar{2}}-\kappa_{1\bar{2}}\wedge\kappa_{2\bar{2}}+\omega_1\wedge\overline{\omega_2}=-(\kappa_{1\bar{1}}-\kappa_{2\bar{2}})\wedge\kappa_{1\bar{2}}.$ This translates into 
   \[(db+(\kappa_{1\bar{1}}-\kappa_{2\bar{2}})b)\wedge dz=0.\]
   The rest follows exactly as in  Lemma \ref{LemmaI}.
 \end{proof}
We say a curve has {\it null-torsion }if $\mathbf{II}=0$. Since $\wedge^2\mathcal{H}\otimes\mathcal{V}\cong \mathbf{C}$ we have 
\[N_X\otimes L\otimes\mathcal{V}\cong \mathbf{C}.\]
If $\mathbf{II}$ is not identically $0$, we define the {\it planar divisor} by
 \[P=\sum_{p: \mathbf{II}(p)=0} ord_p(\mathbf{II})p.\]
In this case, we have 
\[N_X=[P]\otimes L\otimes TM.\]
\begin{theorem}
Let $M=\mathbf{CP}^1$. Then any complex curve $X:M\rightarrow \mathbf{CP}^3$ either is one of the vertical fibers or horizontal or has null-torsion.
\end{theorem}
\begin{proof}

Assume both $\mathbf{I}_1$ and $\mathbf{I}_2$ are not identically $0$. We must show that $\mathbf{II}$ vanishes identically. If not, we have, for $R_1, R_2, P\geq 0$,
\[\mathcal{V}_X=[R_2]\otimes TM, L=[R_1]\otimes TM, N_X=[P]\otimes L\otimes TM,\]
which implies, since $N_X\otimes\L\otimes\mathcal{V}\cong \mathbf{C},$ that 
\[(TM)^3\otimes [2R_1+P+R_2]\cong \mathbf{C},\]
thus ${\text deg   } TM\leq 0$, but ${\text deg   } TM=2$ when $M=\mathbf{CP}^1$.
\end{proof}

\begin{remark}
The computation in this theorem actually shows that if $M^2$ has genus $g$, then any pseudoholomorphic curve $X: M\rightarrow \mathbf{CP}^3$ with none of $\mathbf{I}_1$, $\mathbf{I}_2$ and $\mathbf{II}$ vanishing identically must satisfy 
\[6(g-1)=2 deg(R_1)+deg(R_2)+deg(P).\]
This puts severe restrictions on the bundles $L$, $\mathcal{V}_X$ and $N_X$. For example, if $g=1$, so that $M$ is elliptic, then a pseudoholomorphic curve $X:M\rightarrow \mathbf{CP}^3$ must satisfy $R_1=R_2=P=0$, so that $\mathcal{V}_X=TM$, $L=TM$ and $N_X=(TM)^2$.
\end{remark}

If the pseudoholomorphic curve $X:M^2\rightarrow \mathbf{CP}^3$ has $\mathbf{I}_1\neq 0$, we have a lift of $X$ to a map $\hat{X}:M^2\rightarrow \mathbf{Fl}$ defined by $x\mapsto (X(x), N_X(x)\otimes X(x))$. Some clarification may be necessary. The bundle $N_X$ can be viewed canonically as a subbundle of $X^*(\tilde{Q}\otimes\epsilon^*)\subset \mathbf{C}^4\otimes X^*(\epsilon^*)$. By tensoring with $X^*{\epsilon}$ and canonically identifying $\epsilon\otimes\epsilon^*=\mathbf{C}$ we see that $N_X(x)\otimes X^*(\epsilon)(x)=N_X(x)\otimes X(x)$ is a complex line in $\mathbf{C}^4$. It is easy to see that this line is Hermitian orthogonal to $X(x)\subset \epsilon$ and $X(x)j\subset \epsilon^*$ and thus $\hat{X}$ is well-defined.  Moreover, $X$ has null torsion iff $\hat{X}^*(\kappa_{2\bar{1}})=0$. Composed with $\Pi_2: \mathbf{Fl}\rightarrow \mathbb{CP}^3$, $\hat{X}$ induces a map $Y=\Pi_2\circ\hat{X}: M^2\rightarrow \mathbb{CP}^3$. 
\begin{theorem}
The assignment $X\mapsto Y$ establishes a $1-1$ correspondence between null-torsion pseudoholomorphic curves in $\mathbf{CP}^3$ and nonconstant holomorphic integrals of the holomorphic contact system $\mathbb{V}^*$ on $\mathbb{CP}^3$.
\end{theorem}
\begin{proof}
It is clear from the structure equations that $Y$ is an integral of $\mathbb{V}^*$ if $X$ has null torsion. Conversely, if $Y:M^2\rightarrow \mathbb{CP}^3$ is a nonconstant holomorphic integral of $\mathbb{V}^*$, there exists a unique line bundle $\mathcal{L}\subset\mathbb{H}\subset T\mathbb{CP}^3$ which contains $Y_{*}TM$. We lift $Y$ to a map $\hat{Y}:M^2\rightarrow \mathbf{Fl}$ by $x\mapsto ((\mathbb{H}/\mathcal{L})(x)\otimes Y(x),  Y(x))$. We define the corresponding map $X=\Pi_1\circ\hat{Y}:M\rightarrow \mathbf{CP}^3$. It is clear from the structure equations that such an $X$ has null-torsion. We next show that if we start with null-torsion curve $X:M^2\rightarrow \mathbf{CP}^3$ and run though the procedure $X\rightarrow Y\rightarrow X$ of the above constructions, we arrive at the original curve. In fact the frame adaptations we made before shows we can arrange $\{e_a\}$ so that $\Pi_1^*L_X(x)$ is spanned by $\frac{1}{\sqrt{2}}e_2j\otimes e_1^*$ and $\Pi_1^*N_X(x)$ is spanned by $\frac{1}{\sqrt{2}}e_2\otimes e_1^*$. Thus by definition $Y(x)=[e_2]$. Since $\overline{\omega_1}=0$, $\Pi_2^*\mathcal{L}$ is spanned by $je_1$ from the structure equations. Therefore $\Pi_2^*(\mathbb{H}/\mathcal{L})=[e_1]$ from which we see $\Pi_1(\hat{Y}(x))=X(x)$. We omit the proof that if we start with $Y$ and run the procedure of constructions $Y\rightarrow X\rightarrow Y$ we get $Y$ back.  
\end{proof}
As mentioned before, a powerful construction the integrals of the holomorphic contact system $\mathbb{V}^*$ was provided in \cite{Bryant1} (see Section 3). Of course, there are corresponding results about null-torsion pseudoholomorphic curves in $\mathbf{CP}^3$. We leave most of translation work for the reader and only mention some consequences.
\begin{theorem}
Let $M$ be a compact Riemann surface. There always exists a pseudoholomorphic embedding $M\rightarrow \mathbf{CP}^3$ with null torsion.
\end{theorem}
This is the translation of Theorem G in \cite{Bryant1}.

A horizontal pseudoholomorphic curve $X: M^2\rightarrow \mathbf{CP}^3$ with null torsion corresponds to $Y: M^2\rightarrow\mathbb{CP}^3$ which is {\it superminimal with both positive spin and negative spin} in the sense of Theorem C in \cite{Bryant1}. Thus $M^2$ must be rational. Combining this with Theorem 2.3, we have
\begin{corollary}
There exist pseudoholomorphic curves which are neither vertical nor horizontal.
\end{corollary}
A rational pseudoholomorphic curve is either vertical or horizontal or has null torsion. Both horizontal and null-torsion curves are reduced to integrals of the holomorphic contact system $\mathbb{V}^*$ by Theorem 2.2.  By the result in \cite{Bryant1}, Section 2, such an integral represents a lift of a minimal 2-sphere in $S^4$. Thus the space of nonvertical rational curves in $\mathbb{CP}^3$ can be regarded as the union of 2 copies of the space of minimal 2-spheres in $S^4$. These two copies have a nonempty intersection, corresponding to geodesic 2-spheres.

\end{document}